\documentclass[12pt,a4paper,twoside]{article}
\usepackage[english]{babel}
\usepackage{amssymb}
\usepackage{inputenc}
\newcommand{\ba}{\begin{array}}
\newcommand{\ea}{\end{array}}
\usepackage{amssymb}
\begin{document}

\author{S. Albeverio $^{1},$ Sh. A. Ayupov $^{2,  *},$  K. K.
Kudaybergenov,  $^3$ B. O. Nurjanov $^{4}$}

\title{\bf  Local Derivations on Algebras of Measurable Operators}

\maketitle

\begin{abstract}
The paper is devoted to local derivations on the algebra
$S(\mathcal{M},\tau)$ of $\tau$-measurable operators affiliated
with a von Neumann algebra $\mathcal{M}$ and a faithful normal
semi-finite trace $\tau.$ We prove that every local derivation on
$S(\mathcal{M},\tau)$ which is continuous in the measure topology,
is in fact a derivation. In the particular case of type I von
Neumann algebras they all are inner derivations. It is proved that
for type I finite von Neumann algebras without an abelian direct summand,
and also for von Neumann algebras with the atomic lattice of projections,
the condition of continuity of the local derivation is redundant.
Finally we give necessary and sufficient conditions on a
commutative von Neumann algebra $\mathcal{M}$ for the existence of
local derivations which are not derivations on algebras of
measurable operators affiliated with $\mathcal{M}.$
\end{abstract}

\medskip
$^1$ Institut f\"{u}r Angewandte Mathematik, Universit\"{a}t Bonn,
Wegelerstr. 6, D-53115 Bonn (Germany); SFB 611; HCM;  BiBoS; IZKS; CERFIM
(Locarno); Acc. Arch. (USI), \emph{albeverio@uni-bonn.de}

$^2$ Institute of Mathematics and information  technologies,
Uzbekistan Academy of Sciences, Dormon Yoli str. 29, 100125,
Tashkent (Uzbekistan), e-mail: \emph{sh\_ayupov@mail.ru}

 $^{3}$ Karakalpak state university, Ch. Abdirov str. 1, 142012, Nukus (Uzbekistan),
e-mail: \emph{karim2006@mail.ru}

$^4$ Institute of Mathematics and information  technologies,
Uzbekistan Academy of Sciences, Dormon Yoli str. 29, 100125,
Tashkent (Uzbekistan), e-mail: \emph{nurjanov@list.ru}



* Corresponding author
\newpage

\begin{center}
\textbf{1. Introduction}
\end{center}

The study of derivations on algebras of unbounded operators and in particular on algebras of measurable operators
affiliated with von Neumann algebras is one of the most attractive parts of the general theory of unbounded derivations
on operator algebras.

Given an algebra $A,$ a linear operator $D:A\rightarrow A$ is
called a \textit{derivation}, if $D(xy)=D(x)y+xD(y)$ for all $x,
y\in A$ (the Leibniz rule). Each element $a\in A$ implements a
derivation $D_a$ on $A$ defined as $D_a(x)=ax-xa,\,x\in A.$ Such
derivations $D_a$ are said to be \emph{inner derivations}. If the
element $a,$ implementing the derivation $D_a,$ belongs to a
larger algebra $B$ containing $A,$ then $D_a$ is called \emph{ a
spatial} derivation on $A.$

If the algebra $A$ is commutative, then it is clear that all inner derivations are trivial, i.e. identically zero. One
of the main problem concerning derivations is to prove that a given derivation is inner or spatial, or to show the
existence of non inner (resp. non spatial) derivations and in particular non zero derivations in the commutative cases.

In the paper \cite{Ber} A. F. Ber, V. I. Chilin and F. A. Sukochev
obtained necessary and sufficient conditions for the existence of
non trivial derivations on regular commutative algebras. In
particular it was proved that the algebra $L^{0}(0,1)$ of all
(equivalence classes of) complex measurable function on the $(0,1)$
interval admits non trivial derivations. Independently A. G.
Kusraev \cite{Kus} by the methods of Boolean analysis gave
necessary and sufficient conditions for the existence of non
trivial derivations and automorphisms in extended $f$-algebras. In
particular he also has proved the existence of non trivial
derivations and automorphisms on the algebra $L^{0}(0, 1).$ It is
clear that such derivations are discontinuous and non inner. We
have conjectured in \cite{Alb1},\cite{Alb2} that the existence of
such "pathological"  examples of derivations is closely connected
with the commutative nature of these algebras. This was confirmed
in the particular case of type I von Neumann algebras. Namely in
\cite{Alb1},\cite{Alb2} we have investigated and complety
described derivations on the algebra $LS(\mathcal{M})$ of all
locally measurable operators affiliated with a type I von Neumann
algebra $\mathcal{M}$ and on its various subalgebras. Recently the
above conjecture was also confirmed for the type I case in the
paper \cite{Ber1} by a representation of measurable operators as
operator valued functions. Another approach to similar problems in
$AW ^{*}$-algebras of type I was suggested in the recent paper
\cite{Gut}.

In the paper \cite{Alb3} we have proved have the spatiality of
derivations of the non commutative Arens algebra
$L^{\omega}(\mathcal{M}, \tau)$ associated with an arbitrary von
Neumann algebra $\mathcal{M}$ and a faithful normal semi-finite
trace $\tau.$ Moreover if the trace $\tau$ is finite then every
derivation on $L^{\omega}(\mathcal{M}, \tau)$ is inner.

There exist various types of linear operators which are close to derivations \cite{Bre3},\cite{Bre1},\cite{Kad},\cite{Lar}. In particular R. Kadison \cite{Kad} has introduced and
investigated so called local derivations on von Neumann algebras and some polynomial algebras.

A linear operator $\Delta$ on an algebra $A$ is called a
\textit{local derivation} if given any $x\in A$ there exists a
derivation $D$ (depending on $x$) such that $\Delta(x)=D(x).$  The
main problem concerning this notion is to find conditions under
which local derivations become derivations \cite{Kad} ,
\cite{Lar}. In particular Kadison \cite{Kad} has proved that each
continuous local derivation from a von Neumann algebra
$\mathcal{M}$ into a dual $\mathcal{M}$-bimodule is a derivation.
Later this result was extended in \cite{Bre3} to a larger class of
linear operators $\Delta$ from $\mathcal{M}$ into a normed
$\mathcal{M}$-bimodule $E$ satisfying the identity
\begin{equation}
\Delta(p)=\Delta(p)p + p\Delta(p)
\end{equation}
for every idempotent $p\in \mathcal{M}.$

It is clear that each local derivation satisfies (1) since given any idempotent $p\in \mathcal{M}$ we have
$\Delta(p)=D(p)=D(p^2)=D(p)p + pD(p)= \Delta(p)p + p\Delta(p).$

In \cite{Bre1} it was proved that every linear operator $\Delta$ on the algebra $M_{n}(R)$ satisfying (1) is automatically a derivation,
where $M_{n}(R)$ is the algebra of $n\times n$ matrices over a unital ring $R$ containing $1/2.$

The present paper is devoted to the study of local derivations on the algebra $S(\mathcal{M},\tau)$ of all $\tau$-measurable operators
affiliated with a von Neumann algebra $\mathcal{M}$ and a faithful normal semi-finite trace $\tau.$ The main result (Theorem 2.2) presents
an unbounded version of Kadison's result and it asserts that every local derivation on $S(\mathcal{M},\tau)$ which is continuous in the measure
topology automatically becomes a derivation. In particular in the case of type I von Neumann algebra $\mathcal{M}$ all such local derivations on
$S(\mathcal{M},\tau)$ are inner derivations (Corollary 2.3). We prove also that for type I finite von Neumann algebras without abelian direct summands as well as for von Neumann algebras with the atomic lattice of projections, the
continuity condition on local derivations in Theorem 2.2 is redundant (Theorem 2.5 and Proposition 2.7 respectively).

In section 3 we consider the problem of existence of local
derivations which are not derivations on a class of commutative
regular algebras, which include the algebras of measurable
functions on a finite measure space (Theorem 3.5). As a corollary
we obtain necessary and sufficient conditions for the existence of
local derivation which are not derivations on algebras of
measurable and $\tau$-measurable operators affiliated with a
commutative von Neumann algebra (Theorem 3.8).
 \medskip
\begin{center}
\textbf{2. Continuous local derivations on the algebra $S(\mathcal{M},\tau)$}
\end{center}

Let $H$ be a Hilbert space and let $B(H)$ be the algebra of all bounded linear operators on $H.$
Consider a von Neumann algebra $\mathcal{M}$ in $B(H)$ with a faithful normal semi-finite trace $\tau$. Denote by $P(\mathcal{M})$
the lattice of projections from $\mathcal{M}.$

Recall that a linear subspace $\mathcal{D}$ in $H$ is said to be \textit{affiliated} with $\mathcal{M}$ (and denoted
$\mathcal{D}\eta \mathcal{M}$) if $u(\mathcal{D})\subseteq\mathcal{D}$ for each unitary $u$ from the commutant
$\mathcal{M}'=\{y\in B(H): xy=yx,\,\forall x\in \mathcal{M}\}$ of the algebra $\mathcal{M}.$

A linear operator $x$ acting in $H$ with the domain $\mathcal{D}(x)$ is said to be \textit{affiliated} with $\mathcal{M}$
(denoted $x\eta \mathcal{M}$) if $\mathcal{D}(x) \eta \mathcal{M}$  and $ux(\xi)=xu(\xi)$ for all $u\in \mathcal{M}' $ and $\xi\in\mathcal{D}(x).$

A linear subspace $\mathcal{D}$ in $H$ is said to be \textit{strongly dense} in $H$ with respect to the von Neumann algebra $\mathcal{M},$ if
$\mathcal{D}\eta \mathcal{M}$ and there exists a sequence $\{p_{n}\}^{\infty}_{n=1}$ in $P(\mathcal{M})$ such that $p_n \uparrow \mathbf{1},$
$p_n(H)\subset \mathcal{D}$ and $p_n ^\perp=\mathbf{1}-p_n$ is a finite projection in $\mathcal{M}$ for all $n\in \mathbb{N},$ where $\mathbf{1}$
is the identity in $\mathcal{M}.$

A closed linear operator $x$ acting in the Hilbert space $H$ is said to be \textit{measurable} with respect to $\mathcal{M}$ if
$x\eta \mathcal{M}$ and its domain $\mathcal{D}(x)$ is strongly dense in $H.$

A linear subspace $\mathcal{D}$ in $H$ is called $\tau$-\textit{dense} in $H$ if $\mathcal{D}\eta \mathcal{M}$ and given any $\varepsilon>0$
there exists a projection $p\in \mathcal{M}$ such that $p(H)\subset\mathcal{D}$ and $\tau(p^{\perp})<\varepsilon.$

A linear operator $x$ with the domain $\mathcal{D}(x)\subset H$ is said to be $\tau$-\textit{measurable} with respect to $\mathcal{M}$
if $x\eta \mathcal{M}$ and its domain $\mathcal{D}(x)$ is $\tau$-dense in $H.$

Denote by $S(\mathcal{M})$ and $S(\mathcal{M},\tau)$ respectively the sets of all
measurable and $\tau$-measurable operators affiliated with
$\mathcal{M}$ and consider on $S(\mathcal{M},\tau)$ the topology of convergence in measure (or briefly \textit{measure topology}) $t_{\tau}$
which is given by the following family of neighborhoods of zero:
$$V(\varepsilon,\delta)=\{x\in S(\mathcal{M},\tau):\exists  e\in P(\mathcal{M})\ |\  \tau (e^{\perp})\leq\delta, xe\in
\mathcal{M}, \|xe\|_{\mathcal{M}} \leq\varepsilon \},$$ where $\varepsilon, \delta$ are positive numbers.

It is well-known \cite{Nel} that $S(\mathcal{M},\tau)$ is a complete metrizable topological $*$-algebra with respect to the measure topology $t_{\tau}.$
 \medskip

\textbf{Lemma 2.1.} \textit{The algebra} $S(\mathcal{M},\tau)$ \textit{is semiprime, i.e.}
$aS(\mathcal{M},\tau)a=\{0\}$ \textit{for} $a\in S(\mathcal{M},\tau)$ \textit{implies} $a=0.$

Proof. Let $a\in S(\mathcal{M},\tau)$ and $aS(\mathcal{M},\tau)a=\{0\}$, i.e. $axa=0$ for all
$x\in S(\mathcal{M},\tau).$ In particular for $x=a^{*}$ we have $aa^{*}a=0$ and hence $a^{*}aa^{*}a=0,$ i.e.
$|a|^4=0.$ Therefore $a=0.$ The proof is complete. $\blacksquare$

We are now in position to prove the main result of this section.
 \medskip

\textbf{Theorem 2.2.}  \textit{Let} $\mathcal{M}$ \textit{be a von Neumann algebra with a faithful normal semi-finite trace} $\tau.$
\textit{Then every} $t_{\tau}$-\textit{continuous linear operator} $\Delta$ \textit{on the algebra} $S(\mathcal{M},\tau)$
\textit{satisfying the identity} (1) \textit{is a derivation on} $S(\mathcal{M},\tau).$ \textit{In particular any} $t_{\tau}$-\textit{continuous local
derivation on the algebra $S(\mathcal{M},\tau)$ is a derivation.}

Proof. Given two orthogonal projections $p,\ q \in P(\mathcal{M})$ we have
$$\Delta(p)+\Delta(q)=\Delta((p+q)^2)=\Delta(p+q)(p+q)+(p+q)\Delta(p+q)=$$$$
=\Bigl[\Delta(p)p+p\Delta(p)\Bigr]+\Bigl[\Delta(q)q+q\Delta(q)\Bigr]+\Bigl[\Delta(p)q+p\Delta(q)+\Delta(q)p+q\Delta(p)\Bigr]=$$$$
=\Delta(p)+\Delta(q)+\Bigl[\Delta(p)q+p\Delta(q)+\Delta(q)p+q\Delta(p)\Bigr].$$
Therefore
\begin{equation}
\Delta(p)q+p\Delta(q)+\Delta(q)p+q\Delta(p)=0.
\end{equation}

Denote $ D_{P(\mathcal{M})}=\left\{\sum\limits_{k=1}^{n}\alpha_{k}p_{k}:\alpha_{k}\in \mathbb{R},
 p_{k}\in P(\mathcal{M}),p_{k}p_{l}=0, k\neq l, k,l=\overline{1,n}, n\in \mathbb{N} \right\}$.
For $x=\sum\limits_{i=1}^{n}\alpha_{i}p_{i}\in D_{P(\mathcal{M})}$ we have

$$\Delta(x^2)=
\Delta\left(\left(\sum\limits_{i=1}^{n}\alpha_{i}p_{i}\right)^2\right)=
\Delta\left(\sum\limits_{i=1}^{n}\alpha_{i}^{2}p_{i}\right)=\sum\limits_{i=1}^{n}\alpha_{i}^{2}\Delta(p_{i}),$$ i.e.
\begin{equation}
\Delta(x^2)=\sum\limits_{i=1}^{n}\alpha_{i}^{2}\Delta(p_{i}).
\end{equation}
Further we have
$$\Delta(x)x+x\Delta(x)=
\Delta\left(\sum\limits_{i=1}^{n}\alpha_{i}p_{i}\right)\left(\sum\limits_{i=1}^{n}\alpha_{i}p_{i}\right)+
\left(\sum\limits_{i=1}^{n}\alpha_{i}p_{i}\right)\Delta\left(\sum\limits_{i=1}^{n}\alpha_{i}p_{i}\right)=$$
$$
=\left(\sum\limits_{i=1}^{n}\alpha_{i}^{2}\Delta(p_{i})p_{i}+\sum\limits_{i\neq j}\alpha_{i}\alpha_{j}
\Delta(p_{i})p_{j}\right)+\left(\sum\limits_{i=1}^{n}\alpha_{i}^{2}p_{i}\Delta(p_{i})+
\sum\limits_{i\neq j}\alpha_{i}\alpha_{j}p_{i}\Delta(p_{j})\right)=$$
$$
=\sum\limits_{i=1}^{n}\alpha_{i}^{2}\biggl[\Delta(p_{i})p_{i}+p_{i}\Delta(p_{i})\biggr]+
\sum\limits_{i\neq j}\alpha_{i}\alpha_{j}\biggl[\Delta(p_{i})p_{j}+p_{i}\Delta(p_{j})\biggr]=$$
$$
=\sum\limits_{i=1}^{n}\alpha_{i}^{2}\Delta(p_{i})+
\sum\limits_{i<j}\alpha_{i}\alpha_{j}\biggl[\Delta(p_{i})p_{j}+p_{i}\Delta(p_{j})+\Delta(p_{j})p_{i}+p_{j}\Delta(p_{i})\biggr],
$$
and thus
\begin{equation}
\Delta(x)x+x\Delta(x)=\sum\limits_{i=1}^{n}\alpha_{i}^{2}\Delta(p_{i})+
\sum\limits_{i<j}\alpha_{i}\alpha_{j}\biggl[\Delta(p_{i})p_{j}+p_{i}\Delta(p_{j})+
\Delta(p_{j})p_{i}+p_{j}\Delta(p_{i})\biggr].
\end{equation}
Now (2) and (4) imply
\begin{equation}
\Delta(x)x+x\Delta(x)=\sum\limits_{i=1}^{n}\alpha_{i}^{2}\Delta(p_{i})
\end{equation}
and therefore from (3) and (5) we obtain that $\Delta(x^2)=\Delta(x)x+x\Delta(x)$ for all $x\in D_{P(\mathcal{M})}.$ Since the set
$D_{P(\mathcal{M})}$ is $t_{\tau}$-dense in $S(\mathcal{M},\tau)_{sa}$ and the operator $\Delta$ is $t_{\tau}$-continuous, we have
that  $\Delta(x^2)=\Delta(x)x+x\Delta(x)$ for all $x \in S(\mathcal{M},\tau)_{sa}$ (where $S(\mathcal{M},\tau)_{sa}$ is the space of self-adjoint operators from $S(\mathcal{M},\tau)$).

Now let us show that this relation is valid for arbitrary
operators from $S(\mathcal{M},\tau).$ Consider $x \in
S(\mathcal{M},\tau)$ and let $x=x_{1}+ix_{2}$, where
$x_{1},x_{2}\in S(\mathcal{M},\tau)_{sa}.$ The identity
$$x_{1}x_{2}+x_{2}x_{1}=(x_{1}+x_{2})^2-x^{2}_{1}-x^{2}_{2},$$
implies that
$$\Delta(x_{1}x_{2}+x_{2}x_{1})=\Delta(x_{1})x_{2}+x_{1}\Delta(x_{2})+\Delta(x_{2})x_{1}+x_{2}\Delta(x_{1}).$$
Therefore
$$\Delta(x^2)=\Delta((x_{1}+ix_{2})^2)=\Delta(x_{1}^{2})+i\Delta(x_{1}x_{2}+x_{2}x_{1})-\Delta(x_{2}^{2}),$$ i.e.
$$\Delta(x^2)=\Delta(x_{1}^{2})+i\Delta(x_{1}x_{2}+x_{2}x_{1})-\Delta(x_{2}^{2}).$$
Further we have
$$\Delta(x)x+x\Delta(x)=\Delta(x_{1}+ix_{2})(x_{1}+ix_{2})+(x_{1}+ix_{2})\Delta(x_{1}+ix_{2})=$$
$$=\Bigl[\Delta(x_{1})x_{1}+x_{1}\Delta(x_{1})\Bigr]+i\Bigl[\Delta(x_{1})x_{2}+x_{1}\Delta(x_{2})+
\Delta(x_{2})x_{1}+x_{2}\Delta(x_{1})\Bigr]-
\Bigl[\Delta(x_{2})x_{2}+x_{2}\Delta(x_{2})\Bigr]=$$
$$=\Delta(x_{1}^{2})+i\Delta(x_{1}x_{2}+x_{2}x_{1})-\Delta(x_{2}^{2}),$$
i.e.
$$\Delta(x)x+x\Delta(x)=\Delta(x_{1}^{2})+i\Delta(x_{1}x_{2}+x_{2}x_{1})-\Delta(x_{2}^{2}).$$
Comparing this relation with the above one we obtain
$$\Delta(x^2)=\Delta(x)x+x\Delta(x).$$

This means that $\Delta$ is a Jordan derivation on $S(\mathcal{M},\tau)$ in the sense of \cite{Bre2}. In \cite[Theorem 1]{Bre2} it is proved that
any Jordan derivation on a semiprime algebra is a (associative) derivation. Thus Lemma 2.1 implies that the linear operator $\Delta$ is a
derivation on $S(\mathcal{M},\tau).$ The proof is complete. $\blacksquare$

For type I von Neumann algebras the above result can be strengthened as follows
 \medskip

\textbf{Corollary 2.3.} \textit{Let} $\mathcal{M}$ \textit{be a type I von Neumann algebra with a faithful normal semi-finite trace} $\tau.$
\textit{Then every} $t_{\tau}$-\textit{continuous linear operator} $\Delta$ \textit{on} $S(\mathcal{M},\tau)$
\textit{satisfying} (1) (\textit{in particular every} $t_{\tau}$-\textit{continuous local derivation}) \textit{is an inner derivation.}

Proof. By Theorem 2.2 $\Delta$ is a derivation on
$S(\mathcal{M},\tau).$ By \cite [Corollary 4.5]{Alb2} every
$t_{\tau}$-continuous derivation on $S(\mathcal{M},\tau),$ with
$\mathcal{M}$ of type I, is inner. Therefore $\Delta$ is an inner
derivation on $S(\mathcal{M},\tau).$ The proof is complete.
$\blacksquare$

Now let us show that for finite von Neumann algebras of type I without abelian direct summands the assertion of Corollary 2.3 is valid
for local derivations $\Delta$ without the assumption of $t_{\tau}$-continuity.

Let $\mathcal{M}$ be a homogeneous von Neumann algebra of type I$_{n}, n \in \mathbb{N}$,  with the center $Z$
and with a faithful normal semi-finite trace $\tau.$ In this case $\mathcal{M}$ is $\ast$-isomorphic to the algebra $M_n(Z)$
of all $n\times n$ matrices over $Z,$ and the algebra $S(\mathcal{M}, \tau)$ is $\ast$-isomorphic to the algebra
$M_n(S(Z, \tau_Z))$ of all $n\times n$ matrices over the commutative algebra $S(Z, \tau_Z)$ of $\tau_Z$-measurable operators with respect
to $Z,$ where $\tau_Z$ is the restriction of the trace $\tau$ onto $Z.$ If $\{e_{i,j},\,i,j=\overline{1, n}\}$ is the set of matrix units
in $M_n(S(Z, \tau_Z))$ then each element $x\in M_n(S(Z,\tau_Z))$ is represented as
 $$x=\sum\limits_{i,j=1}^{n}\lambda_{i,j}e_{i,j},\,\lambda_{i,j}\in S(Z, \tau_Z),\,i,j=\overline{1, n}.$$
Let $\delta:S(Z, \tau_Z)\rightarrow S(Z, \tau_Z)$ be a derivation. Setting
\begin{equation}
D_{\delta}(\sum\limits_{i,j=1}^{n}\lambda_{i,j}e_{i,j})=
 \sum\limits_{i,j=1}^{n}\delta(\lambda_{i,j})e_{i,j}
\end{equation}
we obtain a linear operator $D_\delta$ on $M_n(S(Z, \tau_Z)),$ which is a derivation on the algebra
$M_n(S(Z, \tau_Z))$ (see \cite{Alb2}).

Now if $\mathcal{M}$ is an arbitrary finite von Neumann algebra of
type I with the center $Z,$ then there exists a family
$\{z_n\}_{n\in F},$ $F\subseteq\mathbb{N},$ of orthogonal central
projections in $\mathcal{M}$ such that $\sup\limits_{n\in
F}z_{n}=\textbf{1}$ and $\mathcal{M}$ is $\ast$-isomorphic to the
$C^{*}$-product of homogeneous von Neumann algebras $z_n
\mathcal{M}$ of type I$_{n}$ respectively, $n\in F,$ i.e.
$$\mathcal{M}\cong\bigoplus\limits_{n\in F}z_n \mathcal{M}.$$
This implies that the algebra $S(\mathcal{M}, \tau)$ can be embedded as a subalgebra of the direct product of the algebras $S(z_n \mathcal{M}, \tau_n),$
where $\tau_n$ is the restriction of the trace $\tau$ onto the algebra $z_n \mathcal{M},\,n\in F$ (see for details \cite[Section 4]{Alb2}).

Consider a derivation $D$ on the algebra $S(\mathcal{M}, \tau)$ and denote by $\delta$ its restriction on the center $S(Z, \tau_Z)$ of the algebra
$S(\mathcal{M}, \tau).$ Then $\delta$ maps each $z_nS(Z, \tau_Z)\cong Z(S(z_n \mathcal{M}, \tau_n))$ into itself and hence it induces a derivation
$\delta_n$ on $z_nS(Z, \tau_Z)$ for each $n\in F.$

Define as in (6) the derivation $D_{\delta_n}$ on the matrix algebra $M_n(z_nZ(S(\mathcal{M}, \tau)))\cong S(z_n\mathcal{M}, \tau_n)$
for each $n\in F.$ Put
\begin{equation}
D_\delta(\{x_n\}_{n\in F})=\{D_{\delta_n}(x_n)\},\,\{x_n\}_{n\in F}\in S(\mathcal{M}, \tau).
\end{equation}
In \cite{Alb2} it is proved that $D_{\delta}$ is a derivation on the algebra $S(\mathcal{M}, \tau),$ which is restricted to the center of
$S(\mathcal{M}, \tau)$ coincides with $\delta$ (and thus with $D$). In \cite[Lemma 4.3]{Alb2} it has been proved that an arbitrary derivation
$D$ on the algebra $S(\mathcal{M}, \tau)$ for the finite type I von Neumann algebra $\mathcal{M}$ can be uniquely decomposed into the sum
\begin{equation}
D=D_{a}+D_{\delta}
\end{equation}
where $D_{a}$ is an inner derivation on $S(\mathcal{M},\tau)$ implemented by an element $a\in S(\mathcal{M},\tau)$ and $D_{\delta}$ is the
derivation defined as in (7).

Further we shall need the following technical result.
 \medskip

\textbf{Lemma 2.4.} \textit{Every local derivation} $\Delta$ \textit{on the algebra} $S(\mathcal{M},\tau)$
\textit{is necessarily} $P(Z)$-\textit{homogeneous, i.e.} $$\Delta(zx)=z\Delta(x)$$
\textit{for any central projections} $z\in P(Z)=P(\mathcal{M})\cap Z$ and for all $x\in S(\mathcal{M},\tau).$

Proof. Take $z\in P(Z)$ and $x\in S(\mathcal{M},\tau).$ For the element $zx$ by the definition of the local derivation $\Delta$ there exists
a derivation $D$ on $S(\mathcal{M},\tau)$ such that $\Delta(zx)=D(zx).$ Since the projection $z$ is central one has $D(z)=0$ and therefore
$$\Delta(zx)=D(zx)=D(z)x+zD(x)=zD(x),$$ i.e. $\Delta(zx)=zD(x).$ Multiplying by $z$ we obtain $$z\Delta(zx)=z^2D(x)=zD(x)=\Delta(zx)$$ i.e.
$$z^{\bot}\Delta(zx)=(\mathbf{1}-z)\Delta(zx)=0.$$

Therefore by the linearity of $\Delta$ we have $z\Delta(x)=z\Delta(zx)+z\Delta(z^{\bot}x)=z\Delta(zx)=\Delta(zx)$ that is $z\Delta(x)=\Delta(zx).$
The proof is complete. $\blacksquare$
 \medskip

\textbf{Theorem 2.5.} \textit{Let $\mathcal{M}$ be a finite von
Neumann algebra of type I without abelian direct summands and let
$\tau$ be a faithful normal semi-finite trace on $\mathcal{M}.$
Then every local derivation $\Delta$ on the algebra
$S(\mathcal{M},\tau)$ is a derivation and hence can be represented
as} (8).

Proof. Let $\{z_n\}_{n\in F},$ $F\subseteq\mathbb{N}$ be the family of orthogonal central projections in $\mathcal{M}$ with
$\sup\limits_{n\in F}z_{n}=\textbf{1},$ such that $z_n \mathcal{M}$  is a homogeneous von Neumann algebra of type I$_{n}, n\in F.$
Since $\mathcal{M}$ does not have an abelian direct summand we have that $1\not \in F.$

Consider an arbitrary local derivation $\Delta$ on $S(\mathcal{M},\tau).$ By Lemma 2.4 we have that
$$\Delta(z_{n}x)=z_{n}\Delta(x)$$ for each $ n\in F.$ This implies that $\Delta$ maps each
$z_{n}S(\mathcal{M},\tau)=S(z_{n}\mathcal{M},\tau_{n})$ into itself and hence induces a local derivation
$\Delta_{n}=\Delta|_{z_{n}S(\mathcal{M},\tau)}$ on the algebra
$S(z_{n}\mathcal{M},\tau_{n})\cong M_{n}(Z(S(z_{n}\mathcal{M},\tau_{n})))$ for each $n\in F.$ Since $n\neq1,$ \cite [Theorem 2.3]{Bre1}
implies that the operator $\Delta_{n}$ on the matrix algebra $M_{n}(Z(S(z_{n}\mathcal{M},\tau_{n})))$ is a derivation. Therefore
$\Delta=\{\Delta_{n}\}_{n\in F}$ is also a derivation and by \cite[Lemma 4.3]{Alb2} can be uniquely represented in form (8).
The proof is complete. $\blacksquare$
 \medskip

\textbf{Remark 2.6.} In the latter theorem the condition on $\mathcal{M}$ to have no abelian direct summand is crucial, because
in the case of abelian von Neumann algebras the picture is completely different. Local derivations on algebras of $\tau$-measurable
operators affiliated with abelian von Neumann algebras will be considered in the next section.

Now let  $\mathcal{M}$ be a von Neumann algebra with the atomic lattice of projections
and with a faithful normal semi-finite trace
$\tau.$ Then  the von Neumann algebra $\mathcal{M}$ is a direct sum of type
I$_{\alpha}$ factors
$\mathcal{M}_{\alpha}=B(H_{\alpha}),$ where $H_\alpha$ is a Hilbert space
  with $\dim H_\alpha=\alpha$ and the algebra  $S(\mathcal{M}, \tau)$
   can be embedded as a subalgebra of the direct product of the algebras $S(\mathcal{M}_\alpha, \tau_\alpha),$
where $\tau_\alpha$ is the restriction of the trace $\tau$
onto the algebra $\mathcal{M}_\alpha.$ As in Theorem 2.5 by Lemma 2.4 we have that
   every local derivation maps each direct summand
 $S(\mathcal{M}_\alpha, \tau_\alpha)=B(H_{\alpha})$
into itself. Therefore we obtain from \cite[Corollary 3.8]{Bre1} (see also \cite{Lar})
that every local derivation on the
algebra $S(\mathcal{M}, \tau)$ is a derivation and hence by \cite[Corollary 3.11]{Alb1}
and \cite[Theorem 4.4]{Alb2} it is  inner.

Thus, we have proved the following result.

 \medskip

\textbf{Proposition 2.7.} \emph{If M is a von Neumann algebra with the atomic lattice of projections
and with a faithful normal semi-finite trace
$\tau,$ then
every local derivation on the algebra $S(\mathcal{M}, \tau)$ is an inner derivation.}

\begin{center}
\textbf{3. Local derivations on commutative regular algebras}
\end{center}

In this section we shall discuss the problem of existence of local derivations which are not derivations on the algebras $S(\mathcal{M})$
and $S(\mathcal{M},\tau)$ in the case where the von Neumann algebra $\mathcal{M}$ is commutative. Following the approach of the paper
\cite{Ber} we shall consider this problem in a more general setting -- on commutative regular algebras.

Let $A$ be a commutative algebra with the unit $\mathbf{1}$ over
the field $\mathbb{C}$ of complex numbers. We denote by $\nabla$
the set $\{e\in A: e^2=e\}$ of all idempotents in $A.$ For $e,f\in
\nabla$ we set $e\leq f$ if $ef=e.$ With respect to this partial
order, the lattice operation $e\vee f=e+f-ef, \ e\wedge f=ef$ and
the complement $e^{\bot}=\mathbf{1}-e,$ the set $\nabla$ forms a
Boolean algebra. A non zero element $q$ from the Boolean algebra
$\nabla$ is called an \textit{atom} if $0\neq e\leq q, \ e\in
\nabla,$ imply that $e=q.$ If given any nonzero $e\in \nabla$
there exists an atom $q$ such that $q\leq e,$ then the Boolean
algebra $\nabla$ is said to be \textit{atomic}.

An algebra $A$ is called \textit{regular} (in the sense of von Neumann) if for any $a\in A$ there exists $b \in A$ such that $a = aba.$

Further, we shall always assume that $A$ is a unital commutative
regular algebra over $\mathbb{C},$ and that $\nabla$ is the
Boolean algebra of all its idempotents. In this case given any
element $a\in A$ there exists an idempotent $e\in \nabla$  such
that $ea=a,$ and if $ga=a, g\in \nabla$, then $e\leq g.$ This
idempotent is called the \textit{support} of $a$ and denoted by
$s(a).$

Suppose that $\mu$ is a strictly positive countably additive finite measure on the Boolean algebra $\nabla$ of idempotent from $A$
and consider the metric $\rho(a,b)=\mu(s(a-b)),\ a,b\in A.$ From now on we shall  \textit{assume that $(A,\rho)$ is a complete metric space}
(cf. \cite{Ayu1}, \cite{Ber}).
 \medskip

\textbf{Example 3.1.} The most important example of a complete commutative regular algebra $(A,\rho)$ is the algebra
$A=L^0(\Omega)=L^0(\Omega,\Sigma,\mu)$ of all (classes of equivalence of) measurable complex functions on a measure space $(\Omega,\Sigma,\mu),$
where $\mu$ is a finite countably additive measure on $\Sigma,$ and
$\rho(a,b)=\mu(s(a-b))=\mu(\{\omega\in \Omega:a(\omega)\neq b(\omega)\})$
(see for details \cite[Lemma]{Ayu1} and \cite[Example 2.5]{Ber}).
 \medskip

\textbf{Remark 3.2.} If $(\Omega,\Sigma,\mu)$ is a general
localizable measure space, i.e. the (not finite in general)
measure $\mu$ has the finite sum property, then the algebra
$L^0(\Omega,\Sigma,\mu)$ is a unital  regular algebra, but
$\rho(a,b)=\mu(s(a-b))$
 is not a metric in general. But one can
represent $\Omega$ as a union of pair-wise disjoint measurable
sets with  finite measures and thus this algebra is a direct sum
of commutative regular complete metrizable algebras from the above example.

Following \cite{Ber} we call an element $a \in A$ \textit{finitely
valued} (respectively, \textit{countably valued}) if $a =
\sum\limits_{k =1}^n {\alpha _k e_k}$, where $\alpha _k \in
\mathbb{C}$, $e_k \in \nabla, \ e_k e_j=0, \ k\neq j,\ k,j
=1,...,n, \ n\in \mathbb{N}$ (respectively, $a = \sum\limits_{k
=1}^\omega {\alpha _k e_k}$, where $\alpha_k \in \mathbb{C}$, $e_k
\in \nabla, \ e_k e_j=0, \ k\neq j, \ k,j =1,...,\omega,$ where
$\omega$ is a natural number or $\infty$ (in the latter case the
convergence of series is understood with respect to the metric
$\rho$)). We denote by $K(\nabla)$(respectively $K_c(\nabla)$) the
set of all finitely valued (respectively countably valued)
elements in $A.$ It is known that $\nabla \subset K(\nabla)
\subset K_c(\nabla),$ both $K(\nabla)$ and $K_c(\nabla)$ are
regular subalgebras in $A,$ and moreover the closure of
$K(\nabla)$ in $(A,\rho)$ coincides with $K_c(\nabla)$ (see
\cite[Proposition 2.8]{Ber}).

Now let $D$ be a derivation on the given regular commutative algebra $A.$ By \cite[Proposition 2.3]{Ber} (see also \cite[Theorem]{Ayu1})
we have that $s(D(a))\leq s(a)$ for any $a\in A,$ and $D|_{\nabla}=0.$ Therefore by the definition, each local derivation $\Delta$ on
$A$ satisfies the following two condition:
\begin{equation}
s(\Delta(a))\leq s(a), \  \forall a\in A,
\end{equation}
\begin{equation}
\Delta|_{\nabla}\equiv0.
\end{equation}
This means that (9) and (10) are necessary conditions for a linear operator $\Delta$ to be a local derivation on the algebra $A.$
We are going to show that  these two condition are in fact also sufficient.

First we recall some further notions from the paper \cite{Ber}.

Let $B$ be a unital subalgebra in the algebra $A.$ An element $a\in A$ is called:

-- \textit{algebraic with respect to} $B$, if there exists a polynomial $p\in B[x]$ (i.e. a polynomial on $x$ with coefficients
from $B$), such that $p(a)=0$;

-- \textit{integral with respect to} $B$, if there exists a
unitary polynomial $p\in B[x]$ (i.e. the coefficient of the
largest degree of $x$ in $p(x)$ is equal to $\mathbf{1}\in B$),
such that $p(a)=0$;

-- \textit{transcendental with respect to} $B$, if $a$ is not algebraic with respect to $B$;

-- \textit{weakly transcendental with respect to} $B$, if $a\neq 0$ and for any non-zero idempotent $e\leq s(a)$ the element $ea$
is not integral with respect to $B$.
 \medskip

\textbf{Lemma 3.3.} \textit{Given any element $a\in A$  there exists an idempotent $e\in \nabla$  such that}

\textit{$(i)\ ea$ is integral with respect to $K_c(\nabla),$ moreover in this case $ea\in K_c(\nabla)$};

\textit{($ii)\  e^{\bot}a$ is weakly transcendental with respect to $K_c(\nabla),$  if} $e\neq\textbf{1}.$

Proof. Denote by $\nabla_{int}$ the set of all idempotents $e \in \nabla$ such that $ea$ is integral with respect to $K_c(\nabla).$
By \cite[Proposition 3.8]{Ber} each integral element with respect to $K_c(\nabla)$  in fact belongs to $K_c(\nabla).$
Therefore  $\nabla_{int}=\{e\in \nabla:ea\in K_c(\nabla)\}.$ We set $e=\sup\nabla_{int}.$ Since $\nabla$ is a complete Boolean algebra
of countable type \cite[Proposition 2.7]{Ber}, there exists a countable family of mutually disjoint elements $\{e_k\}_{k\geq 1}$ in
$\nabla$ such that $\sup \limits _{k\geq 1}e_k=e$ and given any $e^{\prime}\in \nabla_{int}$ there exists $k\geq 1$ such that
 $e_k\leq e^{\prime}.$ It is clear that $e_k\leq e^{\prime}$,
$e^{\prime}\in \nabla_{int}$ , imply that $e_k\in \nabla_{int}$ and
thus $e_ka\in K_c(\nabla)$. Therefore $ea=\sum \limits _{k\geq
1}e_ka\in K_c(\nabla).$ Further since $s(a)^{\bot}a=0\in
K_c(\nabla)$ we have that $s(a)^{\bot}\leq e,$ i.e. $e^{\bot}\leq
s(a)$ and hence $s(e^{\bot}a)=e^{\bot}.$ Now let us show that if
$e\neq \mathbf{1}$ then $e^{\bot}a$ is weakly transcendental with
respect to $K_c(\nabla).$ Suppose the opposite, i.e. there exists
a non-zero idempotent $q\leq e^{\bot}=s(e^{\bot}a)$ such that $qa$
is integral with respect to $K_c(\nabla).$ This means that $q\in
\nabla_{int},$ i.e. $q\leq e.$  This is a contradiction with
$0\neq q\leq e^{\bot}.$ Therefore $e^{\bot}a$ is weakly
transcendental with respect to $K_c(\nabla)$. The proof is
complete. $\blacksquare$

The following Lemma is the crucial step for the proof of the main results in this section.
 \medskip

\textbf{Lemma 3.4.} \textit{Each linear operator on the algebra $A$ satisfying the conditions} (9) \textit{and} (10)
\textit{is a local derivation on} $A.$

Proof. Let $\Delta$ be a linear operator on the algebra $A$ which satisfies the conditions (9) and (10).
 Let us show that $\Delta|_{K_c(\nabla)}\equiv0.$ Since $\Delta|_{\nabla}\equiv0$ it is clear that $\Delta|_{K(\nabla)}\equiv0.$
Further for $a,b\in A$ we have from (9)
$$\rho(\Delta(a),\Delta(b))=\mu(s(\Delta(a)-\Delta(b)))=
\mu(s(\Delta(a-b)))\leq \mu(s(a-b))=\rho(a,b).$$

This implies that the linear operator $\Delta$ is uniformly continuous with respect to the metric $\rho.$ Since $K(\nabla)$ is dense in
$K_c(\nabla)$ we obtain that $\Delta|_{K_c(\nabla)}\equiv0.$

Now take $a\in A.$ By Lemma 3.3 there exists an idempotent $e\in \nabla$ such that $ea\in K_c(\nabla)$ and $e^{\bot}a$ is weakly
transcendental with respect to $K_c(\nabla).$ Since $\Delta|_{K_c(\nabla)}=0$ and $ea\in K_c(\nabla)$ we have
$$\Delta(a)=\Delta(ea)+\Delta(e^{\bot}a)=\Delta(e^{\bot}a).$$
In particular $s(\Delta(a))=s(\Delta(e^{\bot}a))\leq s(e^{\bot}a).$ Consider the trivial derivation $\delta\equiv 0$ on the regular subalgebra
$K_c(\nabla)$ in $A.$ Now \cite[Proposition 3.7]{Ber} implies that for the weakly transcendental element  $e^{\bot}a$ with respect to the
regular subalgebra $K_c(\nabla)$ and for the element $\Delta(a)$ in $A$ with $s(\Delta(a))\leq s(e^{\bot}a)$ there exists a unique derivation
$\delta_1:B\rightarrow A$ such that $\delta_1(e^{\bot}a)=\Delta(a)$ and $\delta_{1}|_{K_c(\nabla)}\equiv0,$ where $B$ is the subalgebra in $A$
generated by $K_c(\nabla)$ and the element $e^{\bot}a.$ Now by \cite[Theorem 3.1]{Ber} the derivation $\delta_1$ can be extended to a derivation
$D:A\rightarrow A$ and it is clear that $D(e^{\bot}a)=\delta_1(e^{\bot}a)=\Delta(a).$ Further since $ea\in K_c(\nabla)$ and each derivation
satisfies the conditions (9) and (10) we have as above that $D(ea)=0.$ Therefore $$D(a)=D(ea)+D(e^{\bot}a)=D(e^{\bot}a)=\Delta(a),$$ i.e.
for any $a\in A$ we have shown the existence of a derivation $D$ on $A$ such that $D(a)=\Delta(a).$ This means that $\Delta$ is a local derivation
on $A$. The proof is complete. $\blacksquare$

The following is the main result concerning the existence of local derivations on commutative regular algebras.
 \medskip

\textbf{Theorem 3.5. }\textit{Let $A$ be a unital commutative regular algebra over $\mathbb{C}$ and let $\mu$ be a finite strictly
positive countably additive measure on the Boolean algebra $\nabla$ of all idempotents of $A.$ Suppose that $A$ is complete in the
metric $\rho(a,b)=\mu(s(a-b)),\ a,b\in A.$ Then the following conditions are equivalent:}

(i) $K_c(\nabla)\neq A;$

(ii) \textit{the algebra $A$ admits a non-zero derivation;}

(iii) \textit{the algebra $A$ admits a non-zero local derivation;}

(iv) \textit{the algebra  $A$ admits a local derivation which is not a derivation.}

Proof. The implications (i) $\Leftrightarrow$ (ii) are proved in \cite[Theorem 3.2]{Ber}. The assertion (ii) $\Rightarrow$ (iii)
is trivial because any derivation is a local derivation. In order to prove the implication (iii) $\Rightarrow$ (iv) we need the following Lemma.
 \medskip

\textbf{Lemma 3.6.} \textit{If $D$ is a derivation on a commutative regular algebra $A,$ then $D^2$ is a derivation if and only if} $D\equiv 0.$

Proof. Suppose that $D:A\rightarrow A$ is a derivation such that $D^2$ is also a derivation. Then given any $a\in A$ we obtain from the Leibniz
rule for $D^2$ and $D$ respectively:
$$D^2 (a^2)=2aD^2 (a)$$ and $$D^2 (a^2)=D(D(a^2))=D(2aD(a))=2D(a)D(a)+2aD(D(a))=2[D(a)]^2+2aD^2 (a),$$
and therefore $[D(a)]^2=0$.

Since $A$ is regular there exists an element $b\in A$ such that $D(a)=D(a)bD(a).$ Commutativity of $A$ implies that $D(a)=[D(a)]^2b=0,$ i.e.
$D(a)=0$ for all $a\in A.$ The proof of Lemma 3.6 is complete. $\blacksquare$

(iii) $\Rightarrow$ (iv). Since $A$  admits a non-zero local derivation, clearly it admits a non-zero derivation (by the definition of local derivations). From \cite[Theorem 3.2]{Ber} this implies that $K_c(\nabla)\neq A.$ Take an element $a\in A\backslash K_c(\nabla).$
By Lemma 3.3 above there exists an idempotent $e\in \nabla$ such that $ea\in K_c(\nabla)$ and the element $e^{\bot}a$ is weakly transcendental
with respect to $K_c(\nabla)$ provided that  $e\neq \mathbf{1}.$ Since $a\not\in K_c(\nabla),$ we have that $e\neq\mathbf{1},$ and hence the element
$b=e^{\bot}a$ is indeed weakly transcendental with respect to $K_c(\nabla).$ By \cite[Proposition 3.7, Theorem 3.1]{Ber} as in
Lemma 3.4 there exists a derivation $D$ on $A$ such that $D(b)=b.$ Consider the linear map $\Delta=D^2.$ Since $D$ is a derivation on
$A,$ $\Delta$ satisfies the conditions (9) and (10),  and by Lemma 3.4 $\Delta$ is a local derivation on $A,$ and moreover
$\Delta(b)=D(D(b))=D(b)=b,$ i.e. $\Delta \not \equiv 0.$ By Lemma 3.6  $\Delta$ is not a derivation.

(iv) $\Rightarrow$ (i). Let $\Delta$ be a local derivation on $A,$ which is not a derivation. Then it is clear that $\Delta$ is
not identically zero, i.e.  $\Delta(a)\neq 0$ for an appropriate element $a\in A.$ By the definition there exists a derivation
$D$ on $A$ such that $\Delta(a)=D(a)\neq 0,$ i.e. $D$ is a non-zero derivation on $A.$ Therefore by \cite[Theorem 3.2]{Ber}
we obtain that $K_c(\nabla)\neq A.$ The proof of Theorem 3.5 is complete. $\blacksquare$

The important special case of the last theorem is the following result concerning the regular algebra $L^0(\Omega,\Sigma,\mu)$
from the example 3.1.
 \medskip

\textbf{Corollary 3.7.} \textit{Let  $(\Omega,\Sigma,\mu)$ be a finite measure space and let $L^0(\Omega)=L^0(\Omega,\Sigma,\mu)$
be the algebra of all real or complex measurable functions on $(\Omega,\Sigma,\mu).$ The following conditions are equivalent:   }

(i) \textit{the Boolean algebra of all idempotents from $L^0(\Omega)$ is not atomic;}

(ii) $L^0 (\Omega)$ \textit{admits a non-zero derivation;}

(iii) $L^0 (\Omega)$ \textit{admits a non-zero local derivation;}

(iv) $L^0 (\Omega)$ \textit{admits a local derivation which is not a derivation.}

Proof. This follows easily from Theorem 3.5, Example 3.1 and \cite[Theorem 3.3]{Ber}.  $\blacksquare$

It is well known that if $\mathcal{M}$ is a commutative von
Neumann algebra with a faithful normal semifinite trace
$\tau$, then $\mathcal{M}$ is $\ast$-isomorphic to the algebra
$L^{\infty}(\Omega)=L^{\infty} (\Omega,\Sigma,\mu)$ of all essentially
bounded measurable complex function on an appropriate localizable measure space
$(\Omega,\Sigma,\mu)$ and $\tau(f)=\int \limits_{\Omega}f(t)d\mu(t)$
for  $f\in L^{\infty} (\Omega,\Sigma,\mu).$ In this case the algebra $S(\mathcal{M})$ of all
 measurable operators affiliated with
$\mathcal{M}$ may be identified with the algebra $L^0(\Omega)=L^0(\Omega,\Sigma,\mu)$ of all measurable functions on $(\Omega,\Sigma,\mu),$
while the algebra $S(\mathcal{M}, \tau)$ of $\tau$-measurable
operators from $S(\mathcal{M})$ coincides with the algebra
$$\{f\in L^0 (\Omega):\exists F\in\Sigma, \mu(\Omega \setminus F)<+\infty, \chi_F \cdot f \in L^{\infty} (\Omega)\}$$
of all totally $\tau$-measurable functions on $\Omega,$ where $\chi_F$ is the
characteristic function of the set $F.$
If the trace $\tau$ is finite then $S(\mathcal{M})=S(\mathcal{M},\tau)\cong L^0(\Omega)$
 is a commutative regular algebra. But if the trace
$\tau$ is not finite then the algebra $S(\mathcal{M},\tau)$
is not regular.
In this case considering $\Omega$ as a union of pairwise disjoint measurable sets with finite measures we obtain that
$S(\mathcal{M})$ is  a direct sum of commutative regular algebras (see Remark 3.2)
and $S(\mathcal{M},\tau)$ is a subalgebra of the  direct sum of commutative regular algebras.
Therefore from Lemma 2.4 and the above Corollary 3.7 we obtain the following solution of the problem concerning
the existence of non trivial local derivations on algebras of measurable operator in the commutative case.
 \medskip

\textbf{Theorem 3.8.} \textit{Let $\mathcal{M}$ be a commutative von Neumann algebra with a faithful normal
semi-finite trace $\tau.$ The following conditions are equivalent:}

(i) \textit{the lattice $P(\mathcal{M})$ of projections in
$\mathcal{M}$ is not atomic;}

(ii) \textit{the algebra $S(\mathcal{M})$} (\textit{respectively}
$S(\mathcal{M}, \tau)$) \textit{admits a non-inner derivation;}

(iii) \textit{the algebra $S(\mathcal{M})$} (\textit{respectively} $S(\mathcal{M}, \tau)$) \textit{admits
a non-zero local derivation;}

(iv) \textit{the algebra $S(\mathcal{M})$} (\textit{respectively}
$S(\mathcal{M}, \tau)$) \textit{admits a local derivation which is
not a derivation.}

 \medskip

\textbf{Remark 3.9.}  It should be noted that for general (non
commutative) von Neumann algebras the above conditions are not
equivalent but some implications are valid.

The implication (i)
$\Rightarrow$ (ii) is not true in general because for a type
I$_\infty$ von Neumann algebra $\mathcal{M}$ with the non atomic
center $Z$ then the lattice $P(\mathcal{M})$ is not atomic but the algebras $S(\mathcal{M})$ and
$S(\mathcal{M}, \tau)$ do not admit non inner derivations \cite[Lemma 3.5 and Theorem 4.1]{Alb2}.

The implication (ii) $\Rightarrow$ (i) is valid, because if we suppose the lattice $P(\mathcal{M})$
to be atomic then in view of \cite[Corollary 3.1]{Alb1} and \cite[Lemma 3.5 and Theorem 4.1]{Alb2}
every derivation on the algebras $S(\mathcal{M})$ and $S(\mathcal{M}, \tau)$ is automatically $Z$-linear
and hence it is inner, i.e. these algebras do not admit non inner derivations.

The condition (iii) is always fulfilled in the non commutative
case, because every inner derivation which is implemented by a non
central element is a non-zero derivation and hence it is a non-zero
local derivation.

The implications (i) $\Rightarrow$ (iv) and (ii) $\Rightarrow$ (iv) are not true in general, since
if we take a finite von Neumann algebra of type I$ _{n}$  ($n\neq1$),
with a faithful normal finite trace $\tau$ and with  the non atomic center, then
by Theorem 2.5 the algebra $S(\mathcal{M})$ = $S(\mathcal{M}, \tau)$ does not admit a local derivation
which is not a derivation, but it admits non inner derivations of the form $D_{\delta}$ (see Section 2).

The implication (iv) $\Rightarrow$ (i) is true and follows from
Proposition 2.7.

Finally, the implication (iv) $\Rightarrow$ (ii) is an open problem in the general case.

\,
 \newpage

\textbf{Ackowledgments.} \emph{The second and the third named authors would like to acknowledge the hospitality
of the "Institut f$\ddot{u}$r Angewandte Mathematik",
Universit$\ddot{a}$t Bonn (Germany). This work is supported in
part by the DFG 436 USB 113/10/0-1 project (Germany) and the
Fundamental Research Foundation of the Uzbekistan Academy of
Sciences.}

\end{document}